\documentclass{amsart}

\usepackage[nobysame]{amsrefs}
\usepackage{graphicx}
\usepackage{pgf,tikz}
\usetikzlibrary{arrows}

\newcommand{\R}{{\mathbb R}}

\newcommand{\K}{{\mathcal K}}

\newcommand{\conv}{\mathrm{conv}\;}
\newcommand{\vol}{\mathrm{Vol}\;}

\title{The packing density of the $n$-dimensional cross-polytope}

\author{G. Fejes T\'oth}
\address{Alfr\'ed R\'enyi Institute of Mathematics,
Hungarian Academy of Sciences,
Re\'altanoda u. 13-15., 1053, Budapest, Hungary}
\email{gfejes@renyi.hu}

\author{F. Fodor$^\dag$}
\thanks{$\dag$This author was supported by the J\'anos Bolyai Research Scholarship of
the Hungarian Academy of Sciences.}
\address{Department of Geometry, Bolyai Institute, University of Szeged, Aradi v\'ertan\'uk tere 1, 6720 Szeged, Hungary and
Department of Mathematics and Statistics, University of Calgary, Canada}
\email{fodorf@math.u-szeged.hu}

\author{V. V\'{\i}gh$^\ddag$}
\address{Department of Geometry, Bolyai Institute, University of Szeged,
Aradi v\'ertan\'uk tere 1, 6720 Szeged, Hungary}
\email{vigvik@math.u-szeged.hu}
\thanks{$^\ddag$This research was supported by the European Union and the State of Hungary,
co-financed by the European Social Fund in the framework of T\'AMOP-4.2.4.A/ 2-11/1-2012-0001 'National Excellence Program'.}

\date{\today}

\begin{document}

\begin{abstract}
The packing density of the regular cross-polytope in Euclidean $n$-space
is unknown except in dimensions $2$ and $4$ where it is 1. The only
non-trivial upper bound is due to Gravel, Elser, and Kallus \cite{GEK2011}
who proved that for $n=3$ the packing density of the regular
octahedron is at most $1-1.4\ldots\times 10^{-12}$. In this paper,
we prove upper bounds for the packing density of the $n$-dimensional
regular cross-polytope in the case that $n\geq 7$.
We use a modification of Blichfeldt's method \cite{Bl1929}
due to G. Fejes T\'oth and W. Kuperberg \cite{FTGK1993}.
\end{abstract}

\maketitle
\section{Introduction}
Let $K\subset\R^n$ be a convex body (a compact convex set with interior points).
A family $\K=\{K_1, K_2,\ldots\}$ of congruent copies of
$K$ is a {\em packing} in $\R^n$ if the elements of ${\mathcal K}$
are non-overlapping (their interiors are pairwise disjoint).
The {\em density} of a packing is essentially the proportion of
space covered by elements of the packing. The supremum of the
densities of all packings with congruent
copies of a convex body $K$ is called the {\em packing density} of $K$ and
it is denoted by $\delta(K)$. For a more detailed introduction into the basic
properties of density see, for example, \cite{FTL1964} and \cite{Rogers}.

One of the central problems of the theory of packing and covering
is to determine the packing densities of particular convex bodies. The most important
such body has always been the $n$-dimensional unit ball $B^n$. However,
the exact value of $\delta(B^n)$ is known only in the cases $n=2$ and $3$.
In particular, $\delta(B^2)=\pi/\sqrt{12}$, proved by Thue \cites{Thue1892,Thue1910},
and $\delta(B^3)=\pi/\sqrt{18}$, proved by Hales \cite{H2005}. Recently, Cohn and Kumar
\cite{CK2004} proved that in dimension $24$ the density of no sphere packing can exceed the density of the Leech lattice
by more than a multiplicative factor of $1+1.65\times 10^{-30}$. Thus,
for any practical purpose, we may consider $\delta(B^{24})=\pi^{12}/12!$.

The current best asymptotic upper bound for $\delta(B^n)$ is due to Kabatjanski{\u\i} and
Leven\v{s}te{\u\i}n \cite{KL1978}:
\begin{equation}\label{KL}
\delta(B^n)\leq 2^{(-0.599+o(1))n}.
\end{equation}
For small dimensions, the bounds proved by Cohn and Elkies \cite{CE2003}
are better but asymptotically their upper bound is the same as \eqref{KL}.

Other than the $n$-dimensional ball, the most interesting convex bodies are probably the $n$-dimensional
regular polytopes that exist in every dimension: the simplex, the cube, and the cross-polytope. The $n$-cube is a tile, so its packing density is $1$. However, very little is known about the packing densities of the regular simplex and the regular cross-polytope for $n\geq 3$. An exception is the case when $n=4$, then the regular cross-polytope tiles $\R^4$ (cf. Section~22 in \cite{FTL1964}),
thus its packing density is $1$.
Very recently, Gravel, Elser and Kallus
\cite{GEK2011} proved upper bounds for the packing density of the regular tetrahedron
($1-2.6\ldots\times 10^{-25}$) and the regular octahedron ($1-1.4\ldots\times 10^{-12}$).
These bounds are certainly not optimal. It seems unclear whether the method of Gravel, Elser and
Kallus can be extended to higher dimensions.

We note that there has been much work done recently
in order to construct efficient packings of regular tetrahedra, octahedra and other solids in $\R^3$.
For an overview see Torquato and Jiao \cite{TJ2009physrev}.
Dense packings of $n$-dimensional
cross-polytopes were constructed by Rush
\cite{Rush1991}.
Finally, we remark that, to the best of our knowledge, essentially
nothing is known about the covering densities of the regular solids.

In the next section, we will prove an upper bound for
the packing density of the regular $n$-dimensional cross-polytope
using known upper bounds on $\delta (B^n)$ and the ratio of the volumes of the
cross-polytope and its insphere.
Subsequently, we significantly improve these upper bounds
for small dimensions employing a modification of the method of Blichfeldt
by G. Fejes T\'oth and W. Kuperberg \cite{FTGK1993}.
With these methods we establish non-trivial upper bounds on the packing
density of the $n$-dimensional cross-polytope for $n\geq 7$,
and we also show that the packing density of the cross-polytope
approaches $0$ exponentially fast as the dimension tends to infinity.

\subsection{Upper bound using the insphere}
If $K\subset\R^n$ is a convex body and $r(K)$ is the radius of the maximum size
ball contained in $K$, then
\begin{equation}\label{volume-ratio}
\delta (K)\leq \frac{\vol(K)}{r(K)^n \vol (B^n)}\delta (B^n).
\end{equation}
This may provide a non-trivial upper bound on the packing density of
a convex body $K$ whose insphere is sufficiently large in volume
compared to $K$.
Using this idea, Torquato and Jiao \cites{TJ2009, TJ2009physrev} derived
upper bounds for the packing densities of the
regular dodecahedron and the regular icosahedron and
for some of the Archimedean solids in $\R^3$, see Tables III and IV
in \cite{TJ2009physrev}. However, it appears that they
have not used the insphere volume ratio to investigate
the packing densities of convex bodies in higher dimensions.

In $\R^4$, the cube, the cross-polytope and the $24$-cell are tiles,
so their packing densities are all equal to $1$.
The insphere volume ratio method gives an upper bound
greater than $1$ for the
packing density of the regular simplex.
However, one obtains non-trivial
upper bounds for the packing densities of the
$120$-cell and the $600$-cell using the $\delta (B^4)\leq 0.13126\cdot \pi^2/2$
bound by Cohn and Elkies \cite{CE2003}, see the numerical values in
Table~\ref{regular-solids}.

\begin{table}[htb]
$\begin{array}{c|c|c}
n & P  & \text{Upper bound on } \delta(P)\\
\hline
4 & \text{$120$-cell} & 0.74972\\
4 &\text{$600$-cell}  & 0.69073
\end{array}$
\bigskip
\caption{\label{regular-solids}
Upper bounds on the packing densities of $4$-dimensional regular
solids obtained from their insphere volume ratios.
}
\end{table}

In dimensions higher than $4$, there exist only three regular solids, the simplex,
the cube and the cross-polytope. The $n$-dimensional cube is always a tile in $\R^n$,
thus its packing density is $1$. The insphere volume ratio method does not provide a
non-trivial upper bound on the density of the regular simplex in any dimension.
However, it gives an upper bound on the density of the $n$-dimensional
cross-polytope, which approaches $0$ exponentially fast as $n$ tends to infinity.

Consider the regular cross-polytope
$$X^n:=\conv (\pm e_1, \ldots, \pm e_n),$$
where $e_i, i=1, \ldots, n$, are the standard orthonormal basis vectors of $\R^n$.
It is clear that
$$\vol(X^n)=\frac{2^n}{n!},$$
and the inradius of $X^n$ is
$$r_n=r(X^n)=1/\sqrt{n}.$$

We say that two non-negative sequences $f(n)$ and $g(n)$ are asymptotically
equal if $\lim_{n\to\infty}f(n)/g(n)=1$. The asymptotic equality
of $f(n)$ and $g(n)$ will be denoted by $f(n)\sim g(n)$.
We write $f(n)\ll g(n)$ if there exists a positive real
number $c$ such that $f(n)\leq c\cdot g(n)$ for all $n$.

Consider a  packing of cross-polytopes in $\R^n$.
Then \eqref{KL}, \eqref{volume-ratio} and the Stirling formula yield that
\begin{align}\label{vr-bound}
\delta(X^n)&\leq\frac{\vol(X^n)}{r_n^n\vol (B^n)}\delta({B^n})\notag\\
&\leq \frac{\vol(X^n)}{r_n^n\vol (B^n)}2^{-0.599n(1+o(1))}\notag\\
&= \frac{2^n\sqrt{n}^n\Gamma(\frac{n}{2}+1)}{n!\sqrt{\pi}^n}2^{-0.599n(1+o(1))}\notag\\
&\sim \frac{1}{\sqrt{2}}\left (\frac{e}{\pi 2^{0.198}}\right )^{\frac{n}{2}}\\
&\ll 0.86850^n.
\end{align}
Thus, we conclude that
$$\delta(X^n)\to 0, \quad {\text{as} }\quad n\to\infty$$
exponentially fast. Note that the upper bound in \eqref{vr-bound} is asymptotic in nature and it says nothing about the packing density of $X^n$ in specific dimensions.
In order to obtain concrete bounds on $\delta(X^n)$,
we must use specific upper bounds for $\delta(B^n)$.
The Cohn-Elkies bounds on $\delta(B^n)$ for $5\leq n\leq 36$ (Table 3 on page 711 in \cite{CE2003}),
and the almost exact value of $\delta(B^{24})$ by Cohn and Kumar \cite{CK2004}
yield by simple computations the upper bounds for $\delta(X^n)$
shown in Table~\ref{table1}.

\begin{table}[bt]
$\begin{array}{c|c|c|c}
n  & \text{Upper bound on $\delta(X^n)$} & n  & \text{Upper bound on $\delta(X^n)$}\\
\hline
24 & 0.98753 & 31 & 0.67265\\
25 & 0.95416 & 32 & 0.63268\\
26 & 0.90259 & 33 & 0.59472\\
27 & 0.85275 & 34 & 0.55877\\
28 & 0.80476 & 35 & 0.52476\\
29 & 0.75871 & 36 & 0.49264\\
30 & 0.71466
\end{array}$

\bigskip

\caption{Upper bound on $\delta(X^n)$ for $24\leq n\leq 36$ using the
insphere volume ratio.\label{table1}}
\end{table}
We will improve on these bounds in Section~3.

\section{Blichfeldt's method and its extension}
Let $K\subset\R^n$ be a convex body.
A non-negative Lebesgue measurable function
$f:[0,\infty)\to \R$ is a Blichfeldt
gauge for $K$ if it satisfies the following conditions.
\begin{itemize}
\item[i)] $I_n(f):=\int_{\R^n}f(|x|)dx<\infty$.
\item[ii)] If $\{\varphi_i: i=1,2,\ldots \}$
is a set of isometries of $\R^n$ such that the collection
$\{\varphi_i K : i=1,2,\ldots\}$ is a packing,
then for any $x\in\R^n$ it holds that
$\sum_{i=1}^{\infty}f(|\varphi_i^{-1}(x)|)\leq 1$.
\end{itemize}
For technical reasons, we assume the following extra condition on $f$.
\begin{itemize}
 \item[iii)] There exists an $r_0>0$ with $f(r)=0$ for all $r>r_0$.
\end{itemize}
The idea of Blichfeldt \cite{Bl1929} was that if $f$ is a gauge
for a convex body $K$, then
$$\delta(K)\leq\frac{\vol (K)}{I_n(f)}.$$
Blichfeldt applied this idea only to the unit ball in \cite{Bl1929}. He
used the gauge
$$f_0(r)=
\begin{cases}
1-\frac{r^2}{2} &\text{ for } 0\leq r\leq \sqrt{2},\\
0 &\text{ for } r>\sqrt{2}
\end{cases}$$
to show that $\delta(B^n)\le (n+2)2^{-(n+2)/2}$ and noted that a slight
improvement of this bound can be obtained by the gauge
\begin{equation}\label{B-gauge}
f^*(r)=
\begin{cases}
f_0(r) &\text{ for } r\geq 1,\\
1-f_0(2-r) &\text{ for } r\leq 1.
\end{cases}
\end{equation}

For $0\leq\varrho\leq r(K)$
the inner parallel domain of $K$ with radius $\varrho$ is defined as
$$K_{-\varrho}:=\{x\in K: \varrho B^n+x\subseteq K\}.$$
For $x\in\R^n$, let $d(x, K_{-\varrho})$ denote the Euclidean distance
of $x$ from $K_{-\varrho}$.
It is proved in \cite{FTGK1993} that if
$f$ is a Blichfeldt gauge for $B^n$, then for any $0<\rho\leq r(K)$,
$$g_\varrho(x)=f\left (\frac{d(x, K_{-\varrho})}{\varrho}\right )$$
is a Blichfeldt gauge for $K$. Thus, writing
$$G(\varrho)=\int_{\R^n}g_\varrho(x)dx,$$
we have
\begin{equation}\label{FTGK}
\delta(K)\leq \frac{\vol (K)}{G(\varrho)}.
\end{equation}
This method gives
an upper bound on $\delta(K)$
for each $0<\varrho\leq r(K)$.
Our objective is to find, or at least estimate, the best such upper bound.

Let $\kappa_n=\vol (B^n)$.
For a convex body $K\subset\R^n$ and a non-negative
real number $\lambda$, the radius $\lambda$ parallel domain $K_\lambda$
of $K$ is the set of points in $\R^n$ whose distance from $K$
is at most $\lambda$. The Steiner formula for the volume
of $K_\lambda$ can be written in the form
$$\vol (K+\lambda B^n)=\sum_{j=0}^{n} \lambda^{n-j}\kappa_{n-j}V_j(K),$$
where $V_j(K)$, $j=0,\ldots, n$ are the intrinsic volumes of $K$
introduced by McMullen \cite{McM1975}.
Note that $V_n(K)=\vol (K)$ is the volume of $K$,
and $2V_{n-1}(K)=S(K)$ is the surface volume of $K$.

Let $I_0(f):=f(0)$. Using an argument that is very similar to the proof of
Steiner's formula, one obtains that
\begin{equation}\label{polinom}
G(\varrho)=\sum_{j=0}^{n} \varrho^{j}I_j(f)V_{n-j}(K_{-\varrho}).
\end{equation}
If \eqref{polinom} can be calculated or estimated explicitly, then
\eqref{FTGK} provides an upper bound for $\delta(K)$.
We will see in the next section
that it can give better results than the ones we can obtain
from the insphere volume ratio.

The Blichfeldt technique may be used for estimating the packing density of
a convex body for which the intrinsic volumes of its inner
parallel domain can be calculated explicitly
or at least estimated numerically.
Two classes of such bodies were exhibited in
\cite{FTGK1993}: cylinders and the radius $1$ outer parallel domain
of segments. In the next sections, we will apply this method to the
cross-polytope $X^n$.

\section{The case of the $n$-dimensional cross-polytope}
Let $P\subset\R^n$ be a polytope
such that all facets of $P$ are tangent to its insphere. In this case we
say that $P$ is circumscribed around its insphere.
It is not difficult to see that if $P$ is such a polytope, then for
all $0\leq \varrho\leq r(P)$, the radius $\varrho$ inner parallel domain
$P_{-\varrho}$ of $P$ is a polytope that is similar to $P$ with
similarity ratio $(r(P)-\varrho)/r(P)$. Since the $j$th intrinsic volume
is homogeneous of degree $j$, it holds that
$$V_j(P_{-\varrho})=\left (\frac{r(P)-\varrho}{r(P)}\right )^j V_j(P),\quad j=0,\ldots, n.$$
Thus, the right hand side of \eqref{polinom} becomes a polynomial of degree $n$ of
$\varrho$ in the case that $P$ is circumscribed around its insphere, that is,
\begin{equation}\label{G}
G(\varrho)=\sum_{j=0}^{n} \varrho^{j}I_j(f)\left (\frac{r(P)-\varrho}{r(P)}\right )^{n-j} V_{n-j}(P)\quad (0< \varrho\leq (r(P))).
\end{equation}
In particular, $X^n$ is circumscribed about
its insphere. Betke and Henk \cite{Betke-Henk1993} determined the following
formula for the $j$th intrinsic volume of $X^n$:
\begin{equation}\label{intrinsic-volume}
V_j(X^n)=2^{j+1}{n\choose j+1}\cdot \frac{\sqrt{j+1}}{j!}\cdot \gamma (n, j),
\end{equation}
where
\begin{equation*}
\gamma (n, j)=\sqrt{\frac{j+1}{\pi}}\int_{0}^{\infty} e^{-(j+1)x^2}\left (
\frac{2}{\sqrt{\pi}}\int_{0}^{x} e^{-y^2} dy\right )^{d-j-1} dx
\end{equation*}
is the outer angle at a $j$-dimensional face of $X^n$.

Although we only defined $G(\varrho)$ for
$0<\varrho\leq r_n$, it is, in fact, well-defined as a polynomial for all $\varrho$,
and thus its derivatives of all orders exist at $\varrho=0$ and $\varrho=r_n$.
Elementary calculus yields that
\begin{align}\label{Gprime}
G'(r_n)& =-r_n^{n-2}I_{n-1}(f)V_1(X^n)+nr_n^{n-1}I_n(f)\notag\\
& = r_n^{n-2}(nr_nI_n(f)-I_{n-1}(f)V_1(X^n)).
\end{align}

B\"or\"oczky and Henk \cite{BH1999} proved that for any fixed $j$,
\begin{equation*}
\gamma (n,j)\sim \frac{1}{2}\frac{(j+1)!}{\sqrt{j+1}}\frac{(\pi\ln n)^{\frac{j}{2}}}{n^{j+1}},
\text{ as } n\to\infty.
\end{equation*}

In particular,
\begin{equation}\label{V1}
V_1(X^n)\sim \sqrt{\pi}\sqrt{\ln n}, \text{ as } n\to\infty.
\end{equation}

It follows from property iii) of $f$ and the Stirling formula that
\begin{align*}
\frac{I_n(f)}{I_{n-1}(f)} & =\frac{\omega_n}{\omega_{n-1}}\frac{\int_0^\infty f(r)r^{n-1}dr}
{\int_0^\infty f(r)r^{n-2}dr}\\
& \leq \frac{\kappa_n}{\kappa_{n-1}}\frac{n}{n-1}\frac{\int_0^\infty f(r)r^{n-2}r_0 dr}
{\int_0^\infty f(r)r^{n-2}dr}\\
& \sim r_0\sqrt{\frac{2\pi e}{n}} \text{ as $n\to\infty$,}
\end{align*}
and thus
\begin{equation}\label{In}
I_n(f)\ll r_0\sqrt{\frac{2\pi e}{n}} I_{n-1}(f).
\end{equation}

Combining \eqref{Gprime}, \eqref{V1} and \eqref{In}, we obtain that
\begin{align}\label{derivative}
G'(r_n)& \ll r_n^{n-2}I_{n-1}(f)\left (r_0\sqrt{\frac{2\pi e}{n}}nr_n-\sqrt{\pi}\sqrt{\ln n}\right )\\
& = \sqrt{\pi}r_n^{n-2}I_{n-1}(f)\left ( r_0\sqrt{2e}-\sqrt{\ln n}\right )\notag\\
& <0\notag
\end{align}
for sufficiently large $n$. 

Note that $G(r_n)=I_n(f) r_n^n$, and thus by \eqref{FTGK},
\begin{align*}
\delta(X^n) & \leq \frac{\vol (X^n)}{I_n(f) r_n^n}\\
& =\frac{\vol (X^n)}{r_n^n\vol (B^n)}\cdot \frac{\vol (B^n)}{I_n(f)},
\end{align*}
which is exactly the upper bound on $\delta(X^n)$ that we obtain using the
ratio of the volumes of $K$ and its insphere multiplied by the upper
bound on $\delta (B^n)$ from the Blichfeldt gauge $f$.
Thus, \eqref{derivative} implies that for large $n$ the Blichfeldt
method gives a better upper bound on $\delta (X^n)$ than the insphere volume ratio
combined with the upper bound $\delta(B^n)$
coming from the Blichfeldt gauge $f$.

\section{Using the original Blichfeldt gauge function}
In this section, we will use the Blichfeldt gauge $f^*$
as defined in \eqref{B-gauge}. Then $I_0(f^*)=1$, and for $n\geq 1$,
$$I_n(f^*)=\frac{2\kappa_n}{n+2}(\sqrt{2})^n(1+b_n),$$
where
$$b_n=\frac{1}{(\sqrt{2})^n(n+1)}-(\sqrt{2}-1)^{n+1}\left (1+\frac{\sqrt{2}}{n+1}\right ).$$
Routine calculations show that
\begin{align*}
G'(0) & = -\frac{n}{r_n}I_0(f^*)V_n(X^n)+I_1(f^*)V_{n-1}(X^n)\\
& = -n\sqrt{n}\frac{2^n}{n!}+2\sqrt{n}\frac{2^n}{(n-1)!}\frac{1}{2}\\
& = 0,
\end{align*}
 and
\begin{align*}
G''(0) & = \frac{n(n-1)}{r_n^2}I_0(f^*)V_n(X^n)-2\frac{n-1}{r_n}I_1(f^*)V_{n-1}(X^n)
+2I_2(f^*)V_{n-2}(X^n)\\
& = 2I_2(f^*)V_{n-2}(X^n)-\frac{n-1}{r_n}I_1(f^*)V_{n-1}(X^n)\\
& = \frac{n2^n}{(n-2)!}\left (\frac{\sqrt{n-1}}{2}\arccos \left (1-\frac{2}{n}\right )1.062097-1\right ).
\end{align*}
It is easy to check that $G''(0)>0$ when $n\geq 7$. Since $\vol (X^n)/G(0)=1$, our method
provides a non-trivial upper bound on $\delta(X^n)$ in the case that $n\geq 7$.
Furthermore, since $G'(r_n)<0$, the minimum of $\vol (X^n)/G(\varrho)$ is
attained at an interior point of the interval $[0,r_n]$ and the method yields an
upper bound on $\delta(X^n)$ that is
better than the one obtained from the insphere volume ratio combined with the
Blichfeldt upper bound on $\delta(X^n)$.

Although the quantities in \eqref{intrinsic-volume} cannot be
calculated explicitly, they can be approximated by numerical methods.
By such numerical calculations, one obtains for $G(\varrho)$ a degree $n$
polynomial in $\varrho$ whose maximum may be approximated (again
by numerical methods). The results of our calculations are
described in Tables~\ref{table3} and \ref{table4}, and they are
compared to the upper bounds in Table~\ref{table1} in Figure~\ref{figure1}.

\definecolor{qqccqq}{rgb}{0,0.8,0}
\definecolor{ffqqqq}{rgb}{1,0,0}

\begin{figure}[htb]
\centering
\begin{tikzpicture}[line cap=round,line join=round,>=triangle 45,x=.24cm,y=5.0cm]
\draw[->,color=black] (-2.64,0) -- (42.2,0);
\foreach \x in {5,10,15,20,25,30,35,40}
\draw[shift={(\x,0)},color=black] (0pt,2pt) -- (0pt,-2pt) node[below] {\footnotesize $\x$};
\draw[->,color=black] (0,-0.07) -- (0,1.07);
\foreach \y in {,0.2,0.3,0.4,0.5,0.6,0.7,0.8,0.9,1}
\draw[shift={(0,\y)},color=black] (2pt,0pt) -- (-2pt,0pt) node[left] {\footnotesize $\y$};
\draw[color=black] (0pt,-8pt) node[right] {\footnotesize $0$};
\clip(-5.64,-0.07) rectangle (42.2,1.07);
\begin{scriptsize}
\fill [color=ffqqqq] (7,1) circle (1.5pt);
\fill [color=ffqqqq] (8,0.99) circle (1.5pt);
\fill [color=ffqqqq] (9,0.96) circle (1.5pt);
\fill [color=ffqqqq] (10,0.93) circle (1.5pt);
\fill [color=ffqqqq] (11,0.89) circle (1.5pt);
\fill [color=ffqqqq] (12,0.84) circle (1.5pt);
\fill [color=ffqqqq] (13,0.79) circle (1.5pt);
\fill [color=ffqqqq] (14,0.74) circle (1.5pt);
\fill [color=ffqqqq] (15,0.68) circle (1.5pt);
\fill [color=ffqqqq] (16,0.63) circle (1.5pt);
\fill [color=ffqqqq] (17,0.58) circle (1.5pt);
\fill [color=ffqqqq] (18,0.54) circle (1.5pt);
\fill [color=ffqqqq] (19,0.49) circle (1.5pt);
\fill [color=ffqqqq] (20,0.45) circle (1.5pt);
\fill [color=ffqqqq] (21,0.41) circle (1.5pt);
\fill [color=ffqqqq] (22,0.37) circle (1.5pt);
\fill [color=ffqqqq] (23,0.34) circle (1.5pt);
\fill [color=ffqqqq] (24,0.31) circle (1.5pt);
\fill [color=ffqqqq] (25,0.28) circle (1.5pt);
\fill [color=ffqqqq] (26,0.25) circle (1.5pt);
\fill [color=ffqqqq] (27,0.23) circle (1.5pt);
\fill [color=ffqqqq] (28,0.2) circle (1.5pt);
\fill [color=ffqqqq] (29,0.18) circle (1.5pt);
\fill [color=ffqqqq] (30,0.17) circle (1.5pt);
\fill [color=ffqqqq] (31,0.15) circle (1.5pt);
\fill [color=ffqqqq] (32,0.13) circle (1.5pt);
\fill [color=ffqqqq] (33,0.12) circle (1.5pt);
\fill [color=ffqqqq] (34,0.11) circle (1.5pt);
\fill [color=ffqqqq] (35,0.1) circle (1.5pt);
\fill [color=ffqqqq] (36,0.09) circle (1.5pt);
\fill [color=qqccqq] (24,0.99) ++(-1.5pt,0 pt) -- ++(1.5pt,1.5pt)--++(1.5pt,-1.5pt)--++(-1.5pt,-1.5pt)--++(-1.5pt,1.5pt);
\fill [color=qqccqq] (25,0.95) ++(-1.5pt,0 pt) -- ++(1.5pt,1.5pt)--++(1.5pt,-1.5pt)--++(-1.5pt,-1.5pt)--++(-1.5pt,1.5pt);
\fill [color=qqccqq] (26,0.9) ++(-1.5pt,0 pt) -- ++(1.5pt,1.5pt)--++(1.5pt,-1.5pt)--++(-1.5pt,-1.5pt)--++(-1.5pt,1.5pt);
\fill [color=qqccqq] (27,0.85) ++(-1.5pt,0 pt) -- ++(1.5pt,1.5pt)--++(1.5pt,-1.5pt)--++(-1.5pt,-1.5pt)--++(-1.5pt,1.5pt);
\fill [color=qqccqq] (28,0.8) ++(-1.5pt,0 pt) -- ++(1.5pt,1.5pt)--++(1.5pt,-1.5pt)--++(-1.5pt,-1.5pt)--++(-1.5pt,1.5pt);
\fill [color=qqccqq] (29,0.76) ++(-1.5pt,0 pt) -- ++(1.5pt,1.5pt)--++(1.5pt,-1.5pt)--++(-1.5pt,-1.5pt)--++(-1.5pt,1.5pt);
\fill [color=qqccqq] (30,0.71) ++(-1.5pt,0 pt) -- ++(1.5pt,1.5pt)--++(1.5pt,-1.5pt)--++(-1.5pt,-1.5pt)--++(-1.5pt,1.5pt);
\fill [color=qqccqq] (31,0.67) ++(-1.5pt,0 pt) -- ++(1.5pt,1.5pt)--++(1.5pt,-1.5pt)--++(-1.5pt,-1.5pt)--++(-1.5pt,1.5pt);
\fill [color=qqccqq] (32,0.63) ++(-1.5pt,0 pt) -- ++(1.5pt,1.5pt)--++(1.5pt,-1.5pt)--++(-1.5pt,-1.5pt)--++(-1.5pt,1.5pt);
\fill [color=qqccqq] (33,0.59) ++(-1.5pt,0 pt) -- ++(1.5pt,1.5pt)--++(1.5pt,-1.5pt)--++(-1.5pt,-1.5pt)--++(-1.5pt,1.5pt);
\fill [color=qqccqq] (34,0.56) ++(-1.5pt,0 pt) -- ++(1.5pt,1.5pt)--++(1.5pt,-1.5pt)--++(-1.5pt,-1.5pt)--++(-1.5pt,1.5pt);
\fill [color=qqccqq] (35,0.52) ++(-1.5pt,0 pt) -- ++(1.5pt,1.5pt)--++(1.5pt,-1.5pt)--++(-1.5pt,-1.5pt)--++(-1.5pt,1.5pt);
\fill [color=qqccqq] (36,0.49) ++(-1.5pt,0 pt) -- ++(1.5pt,1.5pt)--++(1.5pt,-1.5pt)--++(-1.5pt,-1.5pt)--++(-1.5pt,1.5pt);
\end{scriptsize}
\end{tikzpicture}
\caption{\label{figure1}
Comparison of upper bounds on $\delta (X^n)$ obtained from
insphere volume ratio (diamonds) and the Blichfeldt
method with $f^*$ (dots) for $7\leq n\leq 36$.}
\end{figure}
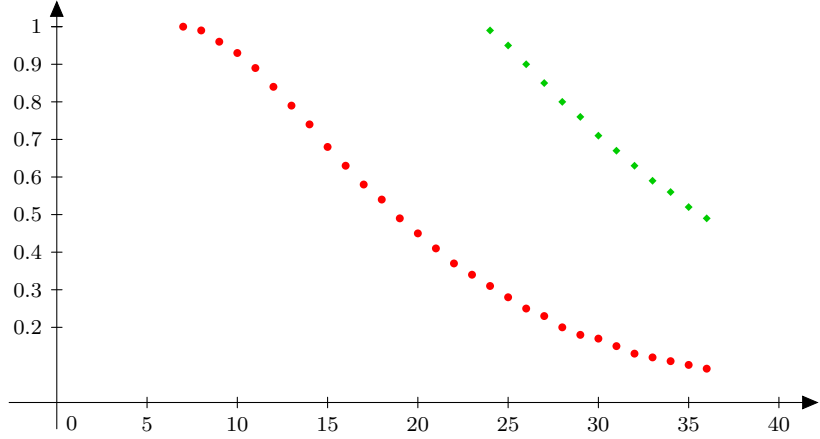

\begin{table}[bth]
$\begin{array}{c|c|c|c}
n & \text{Upper bound} & n & \text{Upper bound} \\
\hline
7 & 0.99805 & 22 & 0.37264\\
8 & 0.98606 & 23 & 0.33850\\
9 & 0.96188 & 24 & 0.30697\\
10 & 0.92730 & 25 & 0.27794\\
11 & 0.88500 & 26 & 0.25129\\
12 & 0.83754 & 27 & 0.22690\\
13 & 0.78705 & 28 & 0.20462\\
14 & 0.73524 & 29 & 0.18448\\
15 & 0.68339 & 30 & 0.16586\\
16 & 0.63247 & 31 & 0.14908\\
17 & 0.58317 & 32 & 0.13398\\
18 & 0.53596 & 33 & 0.12017\\
19 & 0.49116 & 34 & 0.10770\\
20 & 0.44896 & 35 & 0.09647\\
21 & 0.40944 & 36 & 0.08635\\
\end{array}$

\caption{\label{table3} Upper bounds on $\delta(X^n)$ using the
Blichfeldt gauge $f^*$ and
calculated by Maple13.
}
\end{table}

\begin{table}[th]
$\begin{array}{c|c|c|c}
n & \text{density} & n & \text{density} \\
\hline
40 & 5.52108\times 10^{-2} & 140 & 1.98099\times 10^{-7}\\
50 & 1.72421\times 10^{-2} & 150 & 5.36214\times 10^{-8}\\
60 & 5.19017\times 10^{-3} & 160 & 1.44520 \times 10^{-8}\\
70 & 1.52250\times 10^{-3} & 170 & 3.88033 \times 10^{-9}\\
80 & 4.38143\times 10^{-4} & 180 & 1.03837 \times 10^{-9}\\
90 & 1.24242\times 10^{-4} & 190 & 2.77031 \times 10^{-10}\\
100 & 3.48295\times 10^{-5} & 200 & 7.37113 \times 10^{-11}\\
110 & 9.66572\times 10^{-6} & 250 & 9.26781\times 10^{-14}\\
120 & 2.66200\times 10^{-6} & 500 & 2.25312\times 10^{-28}\\
130 & 7.28254\times 10^{-7}& 750 & 4.01494\times 10^{-43}\\
& &  1000 & 6.36493\times 10^{-58}\\

\end{array}$

\caption{\label{table4} Upper bounds on $\delta(X^n)$
using the
Blichfeldt gauge $f^*$ and
calculated by Maple13. 
}
\end{table}

We note that the numerical calculations suggest that
the value of $\varrho$ at which the maximum of $G(\varrho)$ is reached
tends to $\frac{2}{3\sqrt{n}}$ as $n\to\infty$.
Furthermore, it also
appears from calculations that for sufficiently large $n$, the terms in which the exponent of $\varrho$
is around $2n/3$ dominate the polynomial $G(\frac{2}{3\sqrt{n}})$.

Finally, we remark that we fitted and exponential function on the
numerical results obtained from $f^*$ and got the following approximate asymptotics
$$\delta(X^n)\ll 0.87434^n.$$

\section{Concluding remarks}
The upper bounds obtained from the Blichfeldt method
depend on the gauge function used. In Section~4, we used $f^*$ because $I_n(f^*)$
can be evaluated explicitly for all $n$. However, $f^*$ is not the best possible
such function, although in small dimensions it provides better upper bounds on
$\delta (X^n)$  than other known Blichfeldt gauges for $B^n$. For example, it is remarked
in \cite{FTGK1993} that Leven\v{s}te{\u\i}n \cite{L1983} introduced the following Blichfeldt gauge
derived from spherical codes. Let $M(n,\varphi)$ denote the maximum number of
points on $S^{n-1}$ with the property that their pairwise angular distances are
not smaller than $\varphi$. Then the following function
$$
f_n(x)=
\begin{cases}
\frac{1}{M(n,\varphi)} & \text{ for } 0\leq |x|< \sqrt{\frac{2}{1-\cos \varphi}}\\
0 & \text{ for} \sqrt{\frac{2}{1-\cos \varphi}}\leq |x|.
\end{cases}
$$
is a Blichfeldt gauge for $B^n$ if $\pi/3\leq \varphi\leq \pi$.
For a brief explanation why $f_n$ is a Blichfeldt gauge,
see page 726 in \cite{FTGK1993}.
Kabatjansk{\u\i} and Leven\v{s}te{\u\i}n (cf. Formula (52) in \cite{KL1978})
proved that
\begin{equation}\label{KL-precise}
M(n,\varphi)\leq \frac{4{k+n-2\choose k}}{1-t_{1,k+1}^{\alpha, \alpha}}\quad \text{ if }
\cos\varphi\leq t_{1,k}^{\alpha, \alpha},
\end{equation}
where $t_{1,k}^{\alpha, \alpha}$ denotes the largest root of the Jacobi polynomial of degree $k$ with
parameters $\alpha=(n-3)/2$. For a definition of Jacobi polynomials see, for example,
Formula (23) in \cite{KL1978}. Using \eqref{KL-precise}, one can obtain
Blichfeldt gauge functions for $X^n$ in any dimension which yield concrete upper bounds
on $\delta (X^n)$. However, these Blichfeldt gauges do not give better bounds
on $\delta (X^n)$ than $f^*$ up to (at least) dimension $300$. On the other hand,
in dimension $500$ one obtains a better bound using \eqref{KL-precise} than with
$f^*$. We note that the calculations with \eqref{KL-precise} become computationally
very demanding for higher dimensions.

Kabatjansk{\u\i} and Leven\v{s}te{\u\i}n \cite{KL1978} derived 
from \eqref{KL-precise}
\begin{equation}\label{M}
M(n,\varphi)\leq (\sin(\varphi/2))^{-n}2^{-(0.599+o(1))n},
\end{equation}
which holds for $\varphi\leq 63^0$ and is the best asymptotic
upper bound on $M(n,\varphi)$. We note that if one uses
the Blichfeldt gauge $f_n$ that comes from \eqref{M}, then
Blichfeldt's theorem yields the
Kabatjansk{\u\i}--Leven\v{s}te{\u\i}n
upper bound \eqref{KL} on $\delta(B^n)$.
Together with \eqref{derivative}, this indicates that
the Blichfeldt method may provide a better asymptotic upper bound
on $\delta(X^n)$ than the insphere volume ratio combined with \eqref{KL}.

We calculated upper bounds on $\delta(X^n)$ for $n\le1000$ using
the Blichfeldt gauge $f_n$ derived from \eqref{M} (omitting the
unknown $o(n)$ term). The calculations suggest that the minimum of
$\vol (X^n)/G(\varrho)$ is attained at a value $\varrho$ which
tends to roughly $0.767\ldots \times r_n$ as $n\to\infty$.
We fitted an exponential function on the results based
on which we conjecture that
$$\delta(X^n)\ll 0.82886^n.$$
In order to prove this we would need two things: an estimates on the $o(n)$
term in \eqref{M} and asymptotic formulae for
the intrinsic volumes $V_j(X^n)$ for all $j=0,\ldots, n$.
To the best of our knowledge, no such
asymptotic formulae are known at present.


\begin{bibdiv}
\begin{biblist}


\bib{Betke-Henk1993}{article}{
   author={Betke, Ulrich},
   author={Henk, Martin},
   title={Intrinsic volumes and lattice points of cross-polytopes},
   journal={Monatsh. Math.},
   volume={115},
   date={1993},
   number={1-2},
   pages={27--33},
   issn={0026-9255},
   review={\MR{1223242 (94g:52010)}},
}

\bib{Bl1929}{article}{
   author={Blichfeldt, H. F.},
   title={The minimum value of quadratic forms, and the closest packing of
   spheres},
   journal={Math. Ann.},
   volume={101},
   date={1929},
   number={1},
   pages={605--608},
   issn={0025-5831},
   review={\MR{1512555}},
}

\bib{BH1999}{article}{
   author={B{\"o}r{\"o}czky, K{\'a}roly, Jr.},
   author={Henk, Martin},
   title={Random projections of regular polytopes},
   journal={Arch. Math. (Basel)},
   volume={73},
   date={1999},
   number={6},
   pages={465--473},
   issn={0003-889X},
   review={\MR{1725183 (2001b:52004)}},
}

\bib{CE2003}{article}{
   author={Cohn, Henry},
   author={Elkies, Noam},
   title={New upper bounds on sphere packings. I},
   journal={Ann. of Math. (2)},
   volume={157},
   date={2003},
   number={2},
   pages={689--714},
   issn={0003-486X},
   review={\MR{1973059 (2004b:11096)}},
}

\bib{CK2004}{article}{
   author={Cohn, Henry},
   author={Kumar, Abhinav},
   title={The densest lattice in twenty-four dimensions},
   journal={Electron. Res. Announc. Amer. Math. Soc.},
   volume={10},
   date={2004},
   pages={58--67},
   issn={1079-6762},
   review={\MR{2075897 (2005e:11089)}},
}

\bib{CoSl}{book}{
   author={Conway, J. H.},
   author={Sloane, N. J. A.},
   title={Sphere packings, lattices and groups},
   series={Grundlehren der Mathematischen Wissenschaften [Fundamental
   Principles of Mathematical Sciences]},
   volume={290},
   edition={3},
   note={With additional contributions by E. Bannai, R. E. Borcherds, J.
   Leech, S. P. Norton, A. M. Odlyzko, R. A. Parker, L. Queen and B. B.
   Venkov},
   publisher={Springer-Verlag, New York},
   date={1999},
   review={\MR{1662447 (2000b:11077)}},
}

\bib{FTGK1993}{article}{
   author={Fejes T{\'o}th, G.},
   author={Kuperberg, W.},
   title={Blichfeldt's density bound revisited},
   journal={Math. Ann.},
   volume={295},
   date={1993},
   number={4},
   pages={721--727},
   issn={0025-5831},
   review={\MR{1214958 (94c:11055)}},
}

\bib{FTL1964}{book}{
   author={Fejes T{\'o}th, L.},
   title={Regular figures},
   series={A Pergamon Press Book},
   publisher={The Macmillan Co., New York},
   date={1964},
   pages={xi+339},
   review={\MR{0165423 (29 \#2705)}},
}

\bib{GEK2011}{article}{
   author={Gravel, Simon},
   author={Elser, Veit},
   author={Kallus, Yoav},
   title={Upper bound on the packing density of regular tetrahedra and
   octahedra},
   journal={Discrete Comput. Geom.},
   volume={46},
   date={2011},
   number={4},
   pages={799--818},
   issn={0179-5376},
   review={\MR{2846180 (2012g:52035)}},
}

\bib{H2005}{article}{
   author={Hales, Thomas C.},
   title={A proof of the Kepler conjecture},
   journal={Ann. of Math. (2)},
   volume={162},
   date={2005},
   number={3},
   pages={1065--1185},
   review={\MR{2179728 (2006g:52029)}},
}

\bib{KL1978}{article}{
   author={Kabatjanski{\u\i}, G. A.},
   author={Leven{\v{s}}te{\u\i}n, V. I.},
   title={Bounds for packings on the sphere and in space},
   journal={Problems of Information Transmission},
   volume={14},
   date={1978},
   number={1},
   pages={1--17},
   issn={0555-2923},
   review={\MR{0514023 (58 \#24018)}},
}

\bib{L1983}{article}{
   author={Levenshte{\u\i}n, V. I.},
   title={Bounds for packings of metric spaces and some of their
   applications},
   language={Russian},
   journal={Problemy Kibernet.},
   number={40},
   date={1983},
   pages={43--110},
   review={\MR{717357 (86c:52014)}},
}

\bib{McM1975}{article}{
   author={McMullen, P.},
   title={Non-linear angle-sum relations for polyhedral cones and polytopes},
   journal={Math. Proc. Cambridge Philos. Soc.},
   volume={78},
   date={1975},
   number={2},
   pages={247--261},
   issn={0305-0041},
   review={\MR{0394436 (52 \#15238)}},
}

\bib{Rogers}{book}{
   author={Rogers, C. A.},
   title={Packing and covering},
   series={Cambridge Tracts in Mathematics and Mathematical Physics, No. 54},
   publisher={Cambridge University Press, New York},
   date={1964},
   pages={viii+111},
   review={\MR{0172183 (30 \#2405)}},
}

\bib{Rush1991}{article}{
   author={Rush, J. A.},
   title={Constructive packings of cross polytopes},
   journal={Mathematika},
   volume={38},
   date={1991},
   number={2},
   pages={376--380 (1992)},
   issn={0025-5793},
   review={\MR{1147836 (92k:11071)}},
}

\bib{Thue1892}{article}{
  author={Thue, A.},
  title={Om Nogle Geometrisk Taltheoretiske Theoremer},
  journal={Forhdl. Skand. Naturforsk.},
  volume={14},
  pages={352--353},
  date={1892},
}

\bib{Thue1910}{article}{
  author={Thue, A.},
  title={\"Uber die dichteste Zusammenstellung von kongruenten Kreisen in einer Ebene},
  journal={Christiania Vid. Selsk. Skr.},
  volume={1},
  pages={3--9},
  date={1910},
}

\bib{TJ2009}{article}{
    author={Torquato, S.},
    author={Jiao, Y.},
    title={Dense packings of the Platonic and Archimedean solids},
    journal={Nature},
    volume={460},
    date={13 August 2009},
    pages={876--880},
}

\bib{TJ2009physrev}{article}{
   author={Torquato, S.},
   author={Jiao, Y.},
   title={Dense packings of polyhedra: Platonic and Archimedean solids},
   journal={Phys. Rev. E (3)},
   volume={80},
   date={2009},
   number={4},
   pages={041104, 21},
   issn={1539-3755},
   review={\MR{2607418 (2010i:52027)}},
}

\end{biblist}
\end{bibdiv}
\end{document}